\magnification=1200
\def\hal{{\vrule height 10pt width 4pt depth 0pt}}

\def\la{\langle}
\def\ra{\rangle}
\def\tp{{{1}\over{2\pi}}}

\def\C{{\bf C}}
\def\L{{\cal L}}
\def\N{{\bf N}}
\def\R{{\bf R}}
\def\T{{\bf T}}
\def\Z{{\bf Z}}

\centerline{\bf Sub-Riemannian metrics for quantum}

\centerline{\bf Heisenberg manifolds}
\bigskip

\centerline{Nik Weaver}
\bigskip
\bigskip

{\narrower{\it
\noindent Every Heisenberg manifold carries a natural ``sub-Riemannian''
metric with
interesting properties. We describe the corresponding noncommutative metric
structure for Rieffel's quantum Heisenberg manifolds [12].
\bigskip}}
\bigskip

The purpose of this paper is to study an analog, for quantum Heisenberg
manifolds, of the natural sub-Riemannian metrics on classical Heisenberg
manifolds.
\medskip

Quantum Heisenberg manifolds were defined in [12] and they have been further
investigated in [1], [2], [3], and [4]. They are interesting for several
reasons, one being just because they are tractable examples of noncommutative
manifolds. This means that, like the related but simpler noncommutative
tori, quantum Heisenberg manifolds provide a nice concrete setting in which to
explore noncommutative geometry.
\medskip

Our treatment of noncommutative metrics is based on Connes' approach [9].
But we prefer to work with abstract derivations rather than the concretely
presented derivations implicit in Connes' unbounded Fredholm modules. See
[15] for further discussion.
\medskip

Noncommutative metric structure usually arises via an analog of the
classical exterior derivative $d$ on a Riemannian manifold. Classically, this
map may be realized as a derivation from ${\rm Lip}(X) \subset L^\infty(X)$
into the module of bounded measurable 1-forms. The graph of
this derivation is weak*-closed, a property which is characteristic of the
domain being a Lipschitz algebra [15]. In some sense the differentiable
structure resides in the map $d$, while the metric structure resides in
its domain ${\rm Lip}(X)$. There are some examples where one has the latter
sort of structure but not the former ([16], [17]).
\medskip

An interesting feature of the present work is that from the noncommutative
or algebraic point of view, sub-Riemannian metrics are very
close to genuine Riemannian
metrics. In the language of the previous paragraph, the exterior derivative
corresponding to a sub-Riemannian metric is given by composing $d$ with
orthogonal projection onto a closed submodule of $\Omega^1(X)$.
\medskip

This work followed a suggestion by Marc Rieffel, and was helped by
discussions about the Heisenberg group with Daniel Allcock.
\medskip

We adopt the following notational conventions: $c$ is a fixed positive
integer; $\hbar$, $\mu$, and $\nu$ are fixed real numbers; and $H$ is the
Hilbert space $L^2(\R\times\T\times\Z)$.
\bigskip
\bigskip

\noindent {\bf 1. Sub-Riemannian structure for classical Heisenberg manifolds}
\bigskip

Let $M$ be a connected Riemannian manifold. It is well-known that $M$ has a
natural metric such that the distance between two points $x$ and $y$ satisfies
$$d(x,y) = \inf \{l(p): p\hbox{ is a path from $x$ to $y$}\},$$
where $l(p)$ denotes the length of $p$.
\medskip

Now let $B$ be a subbundle of the tangent bundle $TM$. We can use it to define
a new metric $d_B$ by setting
$$\eqalign{d_B(x,y) = \inf\{l(p): p
&\hbox{ is a path from $x$ to $y$}\cr
&\hbox{ which is everywhere tangent to $B$}\}.\cr}$$
This is a {\it sub-Riemannian} or {\it Carnot-Carath\'eodory} metric.
A good general reference on this topic is [6].
Note that we must either require that any two points can be connected by a
path which is tangent to $B$, or else allow distances to be infinite.
\medskip

The simplest non-trivial example of a sub-Riemannian metric arises on the
Heisenberg group. This example is discussed at length in [10].
We now give a brief account of the corresponding construction for Heisenberg
manifolds.
\medskip

The {\it (continuous) Heisenberg group $G$} is the set of all real
$3\times 3$ matrices of the form
$$\pmatrix{1&y&z\cr 0&1&x\cr 0&0&1\cr},$$
with product inherited from the matrix ring $M_3(\R)$. For any positive
integer $c$ the set $H_c$ of elements for which $x$, $y$, and $cz$ are
integers constitutes a discrete subgroup of $G$, and the
quotient construction yields the {\it Heisenberg manifold} $M_c = G/H_c$.
$G$ acts on $M_c$ from the left.
\medskip

$G$ can be identified with $\R^3$ and so it carries a natural differentiable
manifold structure. However, the Euclidean metric on $\R^3$ is not compatible
with the group structure of $G$. Instead, we give
$G$ the unique right-invariant Riemannian metric which agrees with the
Euclidean metric at the origin. Concretely, the three vectors
$${{\partial}\over{\partial x}},\qquad
{{\partial}\over{\partial y}} + x{{\partial}\over{\partial z}},\qquad
{{\partial}\over{\partial z}}$$
define an orthonormal basis at each point $(x,y,z) \in G$.
\medskip

Since this Riemannian metric is right-invariant it descends to $M_c$. The span
of the two vector fields $\partial/\partial x$ and $\partial/\partial y +
x(\partial/\partial z)$ is then a bundle $B$ of tangent planes over $M_c$. (In
fact this is a contact subbundle of $TM_c$, the kernel of the contact 1-form
$\eta = dz - xdy$.) We use it to give $M_c$ a sub-Riemannian metric $d_B$ by
the procedure described above.
\medskip

Interestingly, this metric is finite. That is, even though $M_c$ is
three-dimensional any two points can be joined by a path whose tangent
vector at each point is in the span of $\partial/\partial x$ and
$\partial/\partial y + x(\partial/\partial z)$ [8].
\medskip

Recall that $G$ acts on $M_c$ from the left. The vector fields
$\partial/\partial x$ and $\partial/\partial y + x(\partial/\partial z)$ can
be recovered from this action. To see this
consider the two one-parameter subgroups of
$G$ of the form $x = r$, $y = z = 0$ and $y = s$, $x = z = 0$; the
generators of their actions on $M_c$ are the two desired vector fields.
That is, flowing along the vector fields $\partial/\partial x$ and
$\partial/\partial y + x(\partial/\partial z)$ produces the actions
$\alpha_r$ and $\beta_s$ on $M_c$ defined by
$$\alpha_r(A) = \pmatrix{1&0&0\cr 0&1&r\cr 0&0&1\cr}\cdot A,
\qquad \beta_s(A) = \pmatrix{1&s&0\cr 0&1&0\cr 0&0&1\cr}\cdot A$$
($A \in M_c$). We will use this observation to define an analogous
construction for quantum Heisenberg manifolds.
\bigskip
\bigskip

\noindent {\bf 2. Quantum Heisenberg manifolds}
\bigskip

We recall the definition of the quantum Heisenberg C*-algebras, and define
corresponding von Neumann algebras. Recall that $c$ is a fixed positive
integer, $\hbar$, $\mu$, and $\nu$ are fixed real numbers, and $H =
L^2(\R\times\T\times\Z)$.
\bigskip

\noindent {\bf Definition 1 ([12], Theorem 5.5).}
Let $S^c$ denote the space of $C^\infty$ functions $\Phi$ on $\R\times
\T \times \Z$ which satisfy
\medskip

\noindent (a) $\Phi(x + k, y, p) = e^{ickpy}\Phi(x,y,p)$ for all
$k \in 2\pi\Z$; and
\medskip

\noindent (b) for every polynomial $P$ on $\Z$ and every partial differential
operator $\tilde X = \partial^{m+n}/\partial x^m\partial y^n$ on $\R\times \T$
the function $P(p)\cdot(\tilde X\Phi)(x,y,p)$ is bounded on $C\times \Z$ for
any compact subset $C$ of $\R \times \T$.
\medskip

Define an action of $\Phi \in S^c$ on $H$ by
$$(\Phi\xi)(x,y,p) = \sum_q \Phi(x - \hbar(q - 2p)\mu, y - \hbar(q-2p)\nu, q)
\xi(x, y, p - q)$$
(recall that $H = L^2(\R\times\T\times\Z)$).
Then let $D_\hbar = D_{\hbar, c}^{\mu, \nu}$ be the norm closure of $S^c
\subset B(H)$ and let $N_\hbar =
N_{\hbar, c}^{\mu, \nu}$ be its weak operator closure.\hfill\hal
\bigskip

It is shown in [12] that $D_\hbar$ is a C*-algebra, and it follows that
$N_\hbar$ is a von Neumann algebra. Note that our conventions differ from
[12] by a factor of $2\pi$ in the $x$ variable.
\medskip

The C*-algebras $D_\hbar$ are classified up to isomorphism in [2] and [3].
\medskip

We require alternative characterizations of $D_\hbar$ and $N_\hbar$. The
results are analogous to, but a bit more complicated than,
corresponding facts about noncommutative tori [15].
(Another characterization of $D_\hbar$ is given in [4].)
Our main tool is a kind of Fourier expansion of elements of $N_\hbar$, given
in the next definition. We record its basic properties in the subsequent
lemmas.
\bigskip

\noindent {\bf Definition 2.} For $t \in \R$ and $n \in \Z$ define unitary
operators $U_t$ and $Y_n$ on $H$ by
$$(U_t\xi)(x,y,p) = e^{ipt}\xi(x,y,p)\qquad{\rm and}\qquad
(Y_n\xi)(x,y,p) = \xi(x,y,p+n).$$
For any $T \in B(H)$ and $n \in \Z$ define
$a_n(T) \in B(H)$ by
$$a_n(T) = {{1}\over{2\pi}}\int_{-\pi}^\pi U_tTU_t^{-1}e^{-int}dt.$$
This and all other operator integrals are taken in the weak operator
sense, i.e.\
$$\la a_n(T)\xi, \eta\ra =
{{1}\over{2\pi}}\int_{-\pi}^\pi \la U_tTU_t^{-1} \xi, \eta\ra e^{-int}dt$$
for all $\xi, \eta \in H$.
\medskip

We regard $a_n(T)$ as a sort of Fourier coefficient of $T$; similarly, for
$N \in \N$ we define the Cesaro mean $\sigma_N(T)$ by
$$\sigma_N(T) = \sum_{-N}^N\Bigl(1 - {{|n|}\over{N+1}}\Bigr) a_n(T)
= {{1}\over{2\pi}}\int_{-\pi}^\pi U_tTU_t^{-1} K_N(t)dt,$$
where $K_N$ is the Fej\'er kernel
$$K_N(t) = \sum_{n=-N}^N\Bigl(1-{{|n|}\over{N+1}}\Bigr) e^{-int} =
{{1}\over{N+1}}\Bigl({{\sin((N+1)t/2)}\over{\sin(t/2)}}\Bigr)^2.\eqno{\hal}$$

\noindent {\bf Lemma 3.} {\it For any $T \in B(H)$ we have $\sigma_N(T) \to
T$ weak operator as $N \to \infty$. If the map $t \mapsto U_tTU_t^{-1}$ is
continuous for the norm topology on $B(H)$ then $\sigma_N(T)
\to T$ in norm as $N \to \infty$.}
\medskip

\noindent {\it Proof.} For the weak operator statement, pick $\xi, \eta \in H$
and observe that the map
$$t \mapsto \la (U_tTU_t^{-1} - T)\xi, \eta\ra$$
is continuous and vanishes at $t = 0$. Therefore its integral against
$K_N$ goes to zero as $N \to \infty$ (e.g.\ see [11]), hence
$$\la (\sigma_N(T) - T)\xi, \eta\ra =
\tp\int_{-\pi}^\pi \la (U_t^{-1}TU_t) - T\xi,\eta\ra K_N(t)dt \to 0.$$
This shows that $\sigma_N(T) \to T$ weak operator.
\medskip

For the norm statement, note that the function $t \mapsto
\|U_t^{-1}TU_t^{-1} - T\|$ is continuous and vanishes at zero. Therefore
$$\eqalign{\|\sigma_N(T) - T\|
& = \tp\Big\|\int_{-\pi}^\pi (U_tTU_t^{-1}- T)K_N(t)dt\Big\|\cr
&\leq \tp\int_{-\pi}^\pi \|U_tTU_t^{-1} - T\|K_N(t)dt \to 0\cr}$$
as $N \to \infty$.\hfill\hal
\bigskip

\noindent {\bf Lemma 4.} {\it For any $T \in B(H)$, the operator
$Y_na_n(T)$ preserves constant $p$ subspaces of $H = L^2(\R\times\T\times\Z)$
($p \in \Z$).}
\medskip

\noindent {\it Proof.} Observe that
$Y_na_n(T)$ commutes with $U_s$ for all $s$:
$$\eqalign{Y_na_n(T)U_s
&= \tp Y_n\int_{-\pi}^\pi U_tTU_t^{-1}U_se^{-int}dt\cr
&= \tp Y_nU_s \int_{-\pi}^\pi U_{t - s}TU_{t - s}^{-1}e^{-int}dt\cr
&= \tp e^{ins}U_sY_n\int_{-\pi}^\pi U_tTU_t^{-1}te^{-in(t+s)}dt\cr
&= U_sY_na_n(T).\cr}$$
But the operators $U_s$ generate the von Neumann algebra $l^\infty(\Z) \subset
B(H)$, so we conclude that $Y_na_n(T)$ preserves the constant $p$ subspaces of
$H$.\hfill\hal
\bigskip

\noindent {\bf Lemma 5.} {\it Let $T \in B(H)$.
Suppose $T$ commutes with the operators $V_f$, $W_k$ and $X_r$ defined by
$$\eqalign{(V_f\xi)(x,y,p) &= f(x,y)\xi(x,y,p)\cr
(W_k\xi)(x,y,p) &= e^{-ick(p^2\hbar\nu + py)}\xi(x+k,y,p)\cr
(X_r\xi)(x,y,p) &= \xi(x - 2\hbar r\mu, y - 2\hbar r\nu, p + r)\cr}$$
for all $f \in L^\infty(\R\times \T)$, $k \in 2\pi\Z$, and $r \in \Z$.
Then $a_n(T)$ satisfies
$$(a_n(T)\xi)(x,y,p) = g(x,y,p)\xi(x,y,p-n)$$
for some $g \in L^\infty(\R\times\T\times\Z)$, and the function $g$ satisfies
$$g(x + k, y, p) = e^{-ick((n^2 - 2np)\hbar\nu - ny)}g(x,y,p)\eqno(*)$$
($k \in 2\pi\Z$) and
$$g(x,y,p) = g(x - 2\hbar r\mu, y - 2\hbar r\nu, p + r)\eqno(\dag)$$
($r \in \Z$).
\medskip

The function $\Phi \in L^\infty(\R\times\T\times\Z)$
defined by $\Phi(x,y,p) = 0$ for $p \neq n$ and
$$\Phi(x,y,n) = g(x - \hbar n\mu, y - \hbar n\nu, n)$$
satisfies condition (a) of Definition 1, and $a_n(T) = \Phi$ where $\Phi$
acts on $H$ as in Definition 1.}
\medskip

\noindent {\it Proof.}
Every $V_f$ commutes with $T$ by hypothesis and with both $U_t$ and $Y_n$
by easy computations. It follows that $V_f$ also commutes with $Y_na_n(T)$.
Thus $Y_na_n(T)$, which preserves constant $p$ subspaces by Lemma 4,
must be a multiplication operator on each constant $p$ subspace
of $H$. So $a_n(T)$ has the form
$$(a_n(T)\xi)(x,y,p) = g(x,y,p)\xi(x, y, p - n)$$
for some $g \in L^\infty(\R\times\T\times\Z)$.
\medskip

Next observe that every $W_k$ commutes with both $T$ and $U_t$, hence $W_k$
commutes with $a_n(T)$. So
$$(a_n(T)W_k\xi)(x,y,p) =
g(x,y,p)e^{-ick((p - n)^2\hbar\nu + (p - n)y)}\xi(x + k, y, p - n)$$
equals
$$(W_ka_n(T)\xi)(x,y,p) =
e^{-ick(p^2\hbar\nu + py)}g(x + k, y, p)\xi(x + k, y, p - n),$$
which implies ($*$).
\medskip

Similarly, we have $X_rU_t = e^{irt}U_tX_r$, hence $X_r$ commutes with
$U_tTU_t^{-1}$ and therefore with $a_n(T)$, which implies that
$$(a_n(T)X_r\xi)(x,y,p) =
g(x,y,p)\xi(x - 2\hbar r\mu, y - 2\hbar r\nu, p + r - n)$$
equals
$$(X_ra_n(T)\xi)(x,y,p) =
g(x - 2\hbar r\mu, y - 2\hbar r\nu, p + r)
\xi(x - 2\hbar r\mu, y - 2\hbar r\nu, p + r - n).$$
From this we obtain ($\dag$).
\medskip

Finally, define $\Phi$ as in the statement of the lemma. It follows more or
less immediately from ($*$) that $\Phi$ satisfies condition (a) of
Definition 1. Furthermore it follows from ($\dag$) that the action of $\Phi$
on $H$ given in Definition 1 agrees with the action of $a_n(T)$, that is,
taking $r = n - p$,
$$\eqalign{(a_n(T)\xi)(x,y,p)
&= g(x,y,p)\xi(x,y,p-n)\cr
&= g(x - 2\hbar(n - p)\mu, y - 2\hbar(n - p)\nu, n)\xi(x,y,p - n)\cr
&= \Phi(x - \hbar(n - 2p)\mu, y - \hbar(n - 2p)\nu, n)\xi(x,y,p - n)\cr
&= \sum_q \Phi(x - \hbar(q-2p)\mu, y - \hbar(q-2p)\nu, q)\xi(x,y,p-q).\cr}$$
(Note that even if $\Phi$ does not satisfy condition (b) of Definition 1,
it still acts as a bounded operator
on $H$. In fact $Y_n\Phi$ is a multiplication operator
and so the operator norm of $\Phi$ equals $\|\Phi\|_\infty$.)\hfill\hal
\bigskip

\noindent {\bf Theorem 6.} {\it Let $T \in B(H)$. Then $T \in N_\hbar$ if
and only if $T$ commutes with the operators $V_f$, $W_k$, and $X_r$
defined in Lemma 5 for all $f \in L^\infty(\R\times\T)$, $k \in 2\pi\Z$,
and $r \in \Z$.}
\medskip

\noindent {\it Proof.} The forward direction can be demonstrated by checking
that every $\Phi \in S^c \subset B(H)$ commutes with $V_f$, $W_k$, and $X_r$.
This is an elementary calculation and we omit it.
\medskip

Suppose $T$ commutes with $V_f$, $W_k$, and $X_r$; we must
show that $T \in N_\hbar$. Since $\sigma_N(T) \to T$ weak operator by Lemma
3, it will suffice to show $\sigma_N(T) \in N_\hbar$ for all $N \in \N$;
and since $\sigma_N(T)$ is a linear combination of the Fourier coefficients
$a_n(T)$ it will suffice to show $a_n(T) \in N_\hbar$ for all $n \in \N$.
\medskip

Let $g$ and $\Phi$ be the functions defined in Lemma 5 and
let $(h_m)$ be a sequence of functions in $C^\infty(\R \times\T)$ with
${\rm supp}(h_m) \subset [-1/m, 1/m] \times [-1/m, 1/m]$ and $\|h_m\|_1 = 1$.
Define smoothings $\Phi_m$ of $\Phi$ by twisted convolution:
$$\Phi_m(x,y,n) =
\int_{\R\times\T} h_m(r,s)\Phi(x-r,y-s,n)e^{icnxs}\,drds.$$
Then $\Phi_m$ satisfies condition (a) of Definition 1 because
$$\eqalign{\Phi_m(x + k, y, n)
& = \int h_m(r,s)\Phi(x + k - r, y-s, n)e^{icn(x+k)s}\,drds\cr
& = \int h_m(r,s)e^{ickn(y-s)}\Phi(x - r, y - s, n)e^{icn(x+k)s}\,drds\cr
& = e^{ickny}\Phi_m(x, y, n);\cr}$$
and $\Phi_m$ satisfies condition (b) of Definition 1 because it is $C^\infty$
and supported on $p = n$. So $\Phi_m \in S^c$, and a standard
application of Fubini's theorem and dominated convergence shows that
$\Phi_m \to \Phi$ weak* in
$L^\infty(\R\times\T\times\Z)$. Thus the multiplication operators
$Y_n\Phi_m$ converge weak operator to $Y_n\Phi$, so that $\Phi_m \to \Phi$
weak operator. This shows that $a_n(T) = \Phi$ belongs to $N_\hbar$.\hfill\hal
\bigskip

The corresponding characterization of the C*-algebra $D_\hbar$ can be stated
most naturally in terms of the action of the Heisenberg group $G$ on
the $D_\hbar$ given in [12].
\bigskip

\noindent {\bf Definition 7 ([12], p.\ 557).} For $r, s, t \in \R$ let
$U_{(r,s,t)}$ be the unitary operator on $H$ defined by
$$(U_{(r,s,t)}\xi)(x,y,p) =
e^{ip(t + cs(x + \hbar p\mu - r))}\xi(x - r, y - s, p).$$
Then $L_{(r,s,t)}(T) = U_{(r,s,t)}TU_{(r,s,t)}^{-1}$
defines an action $L$ of $G$ on $B(H)$, where
we take $(r, s, t) \in \R^3 \cong G$. A short computation shows that this
action preserves $S^c$, hence it preserves $D_\hbar$ and $N_\hbar$.
\medskip

By specializing to the three coordinate axes in $G \cong \R^3$
we get three actions $\alpha$, $\beta$, $\gamma$ of $\R$, defined by
$${\eqalign{\alpha_r(T) &= U_{(r,0,0)}TU_{(r,0,0)}^*\cr
\beta_s(T) &= U_{(0,s,0)}TU_{(0,s,0)}^*\cr
\gamma_t(T) &= U_{(0,0,t)}TU_{(0,0,t)}^*.\cr}}$$
(Note that $U_{(0,0,t)}$ equals the unitary $U_t$ of Definition 2.
Thus, for instance,
$a_n(T) = \tp\int_{-\pi}^\pi\gamma_t(T)e^{-int}dt$.)\hfill\hal
\bigskip

For $\Phi \in S^c \subset B(H)$ we have
$$\eqalign{\alpha_r(\Phi)(x,y,p) &= \Phi(x - r, y, p)\cr
\beta_s(\Phi)(x,y,p) &= e^{ipscx}\Phi(x, y - s, p)\cr
\gamma_t(\Phi)(x,y,p) &= e^{ipt}\Phi(x,y,p).\cr}$$

\noindent {\bf Lemma 8} {\it Let $\Phi \in S^c$. Then there exists
$K > 0$ such that
$$\|\alpha_r(\Phi) - \Phi\| \leq Kr\qquad{\rm and}\qquad
\|\beta_s(\Phi) - \Phi\| \leq Ks$$
for all $r, s > 0$.}
\medskip

\noindent {\it Proof.} By condition (b) of Definition 1 there exists a
positive function $f \in L^1(\Z)$ such that
$$\Big|{{\partial\Phi}\over{\partial x}}(x,y,p)\Big| \leq f(p)$$
for all $p \in \Z$. So for $r > 0$ we have
$$|(\alpha_r(\Phi) - \Phi)(x,y,p)| = |\Phi(x - r, y, p) - \Phi(x, y, p)|
\leq r\cdot f(p).$$
It follows that $\|\alpha_r(\Phi) - \Phi\| \leq r\cdot \|f\|_1$.
\medskip

Now use condition (b) of Definition 1 to find a positive function $f_1$
such that $pf_1(p) \in L^1(\Z)$ and $|\Phi(x,y,p)| \leq f_1(p)$ for all
$p \in \Z$ and a positive function $f_2 \in L^1(\Z)$ such that
$$\Big|{{\partial\Phi}\over{\partial y}}(x,y,p)\Big| \leq f_2(p)$$
for all $p \in \Z$. Also note that
$$\sup_{0 \leq x < 2\pi} |e^{ipcsx} - 1| \leq 2\pi cps.$$
Then
$$\eqalign{|(\beta_s(\Phi) - \Phi)(x,y,p)|
&= |e^{ipcsx}\Phi(x, y - s, p) - \Phi(x,y,p)|\cr
&\leq |e^{ipcsx} - 1|f_1(p) + |\Phi(x,y-s,p) - \Phi(x, y, p)|\cr
&\leq (2\pi cpf_1(p) + f_2(p))s\cr}$$
for all $x \in [0,2\pi)$. But for any $x \in [0,2\pi)$ and $k \in 2\pi\Z$
we have
$$\eqalign{(\beta_s(\Phi) - \Phi)(x + k, y, p)
&= e^{ipcs(x+k)}\Phi(x + k, y - s, p) - \Phi(x + k, y, p)\cr
&= e^{ipcs(x+k)}e^{ickp(y-s)}\Phi(x, y - s, p) - e^{ickpy}\Phi(x, y, p)\cr
&= e^{icpky}(\beta_s(\Phi) - \Phi)(x,y,p),\cr}$$
so that the previous estimate holds for all $x \in \R$. It follows that
$\|\beta_s(\Phi) - \Phi\| \leq (2\pi c\|pf_1\|_1 + \|f_2\|_1)s$.\hfill\hal
\bigskip

\noindent {\bf Theorem 9.} {\it
Let $T \in N_\hbar$. Then $T \in D_\hbar$ if and only if
the maps $r \mapsto \alpha_r(T)$ and $s \mapsto \beta_s(T)$ are continuous
for the norm topology on $N_\hbar$.}
\medskip

\noindent {\it Proof.} ($\Rightarrow$) The set of operators $T$ for which
$\alpha$ and $\beta$ are norm-continuous is easily seen to be norm-closed.
Thus it suffices to show that every $\Phi \in S^c \subset
B(H)$ has this property. This was shown in Lemma 8.
\medskip

($\Leftarrow$) Let $T \in N_\hbar$ and suppose $\alpha$ and $\beta$ are
norm-continuous for $T$. It follows that $\gamma$ is also
norm-continuous for $T$ by the identity
$$\gamma_t = \beta_{t'}^{-1}\alpha_{t'}^{-1}
\beta_{t'}\alpha_{t'}$$
where $t' = \sqrt{t/c}$. Thus $\sigma_N(T) \to T$ in norm as $N \to \infty$
by Lemma 3. Therefore, to prove that $T \in D_\hbar$ it will suffice
to show that $\sigma_N(T) \in D_\hbar$, or indeed that $a_n(T) \in D_\hbar$.
\medskip

Now $\alpha$ and $\gamma$ commute, so
$${\eqalign{\|\alpha_r(a_n(T)) - a_n(T)\|
&= \tp\| \int_{-\pi}^\pi (\alpha_r(\gamma_t(T)) - \gamma_t(T))e^{-int}dt\|\cr
& \leq \tp\int_{-\pi}^\pi \|\alpha_r(\gamma_t(T)) - \gamma_t(T)\| dt\cr
& = \tp\int_{-\pi}^\pi \|\alpha_r(T) - T\| dt\cr
& = \|\alpha_r(T) - T\|.\cr}}$$
Since $T$ is norm-continuous for $\alpha$, this shows that $a_n(T)$ is as
well; the same argument shows that $a_n(T)$ is norm-continuous for $\beta$.
\medskip

Using the fact that $\alpha_r(Y_na_n(T)) = Y_n\alpha_r(a_n(T))$
we get that $\alpha_r(Y_na_n(T))$ is continuous in norm as a function
of $r$. Similarly, a short
computation shows that $\beta_s(Y_n)$ is continuous in norm as a function
of $s$, hence $\beta_s(Y_na_n(T)) = \beta_s(Y_n)\beta_s(a_n(T))$
is also continuous in norm. It is then standard that
the functions $g$ and $\Phi$ defined in Lemma 5 must be
uniformly continuous. We let $\Phi$ act on $H$ by the formula given in
Definition 1, so that $a_n(T) = \Phi$ as operators.
\medskip

Now, just as in the proof of Theorem 6, we can smooth $\Phi$ by taking a
twisted convolution with a $C^\infty$ approximate unit of $L^1$-norm one,
to get a sequence $(\Phi_m)$ in $S^c$. But since $\Phi$ is continuous,
$\Phi_m \to \Phi$ in sup norm = operator norm,
hence $a_n(T) = \Phi \in D_\hbar$.\hfill\hal
\bigskip

\noindent {\bf Corollary 10.} {$D_\hbar$ consists of precisely the elements
of $N_\hbar$ for which the action of $G$ is norm-continuous.}
\medskip

\noindent {\it Proof.} Norm-continuity for $G$ implies norm-continuity for
$\alpha$ and $\beta$, so one direction follows immediately from Theorem 9.
For the other direction, we also know from Theorem 9 that every element of
$D_\hbar$ is norm-continuous for $\alpha$ and $\beta$, and that this implies
norm-continuity for $\gamma$ as well. But $\alpha$, $\beta$, and $\gamma$
generate $G$, so this is enough.\hfill\hal.
\bigskip
\bigskip

\noindent {\bf 3. The noncommutative sub-Riemannian metric}
\bigskip

There is a natural noncommutative sub-Riemannian metric on $N_\hbar$, and it
can be presented
in both local and global forms. The local version is a sort of noncommutative
exterior derivative, while the global version is the noncommutative Lipschitz
algebra that is the former's domain.
\bigskip

\noindent {\bf Definition 11.} Let $E = N_\hbar \oplus N_\hbar$; we regard it
as a Hilbert bimodule over $N_\hbar$ with left and right $N_\hbar$-valued
inner products given by
$$\la x_1\oplus x_2, y_1 \oplus y_2\ra_l = x_1y_1^* + x_2y_2^*$$
and
$$\la x_1\oplus x_2, y_1 \oplus y_2\ra_r = x_1^*y_1 + x_2^*y_2.$$
Let $\delta_1$ and $\delta_2$ be the generators of the actions $\alpha$
and $\beta$ defined in the last section, i.e.\
$$\delta_1(x) = \lim_{r \to 0} {{\alpha_r(x) - x}\over{r}}\qquad{\rm and}
\qquad\delta_2(x) = \lim_{s \to 0} {{\beta_s(x) - x}\over{s}}$$
for all $x \in N_\hbar$
for which the limits exist in the weak operator sense. Define $\L_\hbar =
{\rm dom}(\delta_1) \cap {\rm dom}(\delta_2)$ and define $d: \L_\hbar \to E$
by $d(x) = \delta_1(x) \oplus \delta_2(x)$. Give $\L_\hbar$ the norm
$$\|x\|_L = {\rm max}(\|x\|, \|d(x)\|_l, \|d(x)\|_r)$$
where $\|\cdot\|_l$ and $\|\cdot\|_r$ are the left and right Hilbert module
norms on $E$.\hfill\hal
\bigskip

In the case $\hbar = 0$, $\delta_1$ and $\delta_2$ are genuine partial
derivatives, and $d(x)$ is the projection of the exterior derivative of $x$
onto the cotangent subbundle dual to $B$.
\medskip

The following alternative characterization of $\L_\hbar$ is useful. It
follows immediately from ([7], Proposition 3.1.6).
\bigskip

\noindent {\bf Lemma 12.} {\it Let $x \in N_\hbar$. Then $x \in \L_\hbar$
if and only if $\sup_{r > 0}\|\alpha_r(x) - x\|/r$ and
$\sup_{s > 0}\|\beta_s(x) - x\|/s$ are finite.\hfill\hal}
\bigskip

\noindent {\bf Theorem 13.} {\it The map $d$ and its domain $\L_\hbar$ have
the following properties:
\medskip

{\narrower{
\noindent (a). $d$ is an unbounded derivation with weak*-closed graph.
\medskip

\noindent (b). $\L_\hbar$ is a dual Banach algebra. It contains $S^c$ and is
densely contained in $D_\hbar$.
\medskip

\noindent (c). If $\hbar = 0$ then $\L_\hbar$ is naturally identified with
the algebra of functions on $M_c$ which are Lipschitz for the sub-Riemannian
metric defined in section 1, and $\|d(x)\|_l = \|d(x)\|_r$ equals the
Lipschitz number of $x$, for any $x \in \L_0 \cong {\rm Lip}(M_c)$.
\medskip}}}

\noindent {\it Proof.} (a). The fact that $d$ is a derivation, i.e.\ is
linear and satisfies the Liebnitz formula $d(xy) = xd(y) + d(x)y$, is an
elementary calculation. Weak*-closure of the graph follows from ([7],
Proposition 3.1.6).
\medskip

(b). The norm $\|\cdot\|_L$ equals the graph norm when $\L_\hbar$ is identified
with the graph of $d$ by the map $x \mapsto x \oplus d(x)$ and $E$ is given
the max of its left and right Hilbert norms (which is equivalent to the von
Neumann algebra norm on $E$). Thus $\L_\hbar$ is isometric to a
weak*-closed subspace of a dual Banach space, hence $\L_\hbar$ is a dual space.
It is an algebra because ${\rm dom}(\delta_1)$ and ${\rm dom}(\delta_2)$
are (being domains of derivations).
\medskip

$S^c \subset \L_\hbar$ follows from Lemma 8 and Lemma 12.
$\L_\hbar \subset D_\hbar$ follows from Theorem 9 and Lemma 12, using the
fact that continuity of an $\R$-action is equivalent to continuity at 0. Also,
$\L_\hbar$ is dense in $D_\hbar$ because it contains $S^c$.
\medskip

(c). Let $\hbar = 0$. It is straightforward to check that $D_0$ and
$N_0$ are, respectively, naturally isomorphic to $C(M_c)$ and
$L^\infty(M_c)$; this simply involves taking the Fourier transform in the
$p$ variable. Now if $f \in L^\infty(M_c)$ is Lipschitz
for the sub-Riemannian metric then it satisfies
$$\|f\circ\phi - f\|_\infty \leq L(f)\cdot r$$
for any isometry $\phi$ of $M_c$ such that $d_B(\phi(\rho), \rho) = r$ for all
$\rho \in M_c$, where $L(f)$ is the Lipschitz number of $f$. Taking
$$\phi(r,s,t) = (r - h\cos\theta, s - h\sin\theta, t - hr\sin\theta),$$
this shows that ${\rm Lip}(M_c) \subset \L_0$ by Lemma 12; and
$$\|d(f)\|_l(\rho) = {\rm ess}\sup_\theta |\delta_1(f)(p)|^2\cos^2\theta +
|\delta_2(f)(p)|^2\sin^2\theta \leq L(f)$$
for almost every $\rho \in M_c$, so $\|d(f)\| \leq L(f)$.
\medskip

Conversely, let $f$ be any function in $\L_0$ and let $\rho,\sigma \in M_c$.
By ([5], Theorem 2.7) there exists a constant velocity geodesic
$p: [0,1] \to M_c$
which is everywhere tangent to the subbundle $B$ of $TM_c$ defined in
section 1 and which satisfies $p(0) = \rho$, $p(1) = \sigma$, and
$l(p) = d_B(\rho,\sigma)$. Then the
function $g = f\circ p: [0,1] \to \C$ satisfies $g(0) = f(\rho)$, $g(1) =
f(\sigma)$, and $g'(t) = \la d(f)(t), dp(t)\ra$ for almost every $t \in [0,1]$.
This implies that $g$ is Lipschitz with
$$L(g) \leq \|dp\| = \|d(f)\|\cdot l(p),$$
so
$${{|f(\rho) - f(\sigma)|}\over{d_B(\rho,\sigma)}} =
{{g(0) - g(1)}\over{l(p)}} \leq L(g)/l(p) \leq \|d(f)\|.$$
Taking the supremum over all $\rho$ and $\sigma$
shows that $f$ is Lipschitz and $L(f) \leq \|d(f)\|$.\hfill\hal
\bigskip

We conclude this section with a proof that the unit ball of $\L_\hbar$ is
compact in operator norm. In the commutative case, this is true of
${\rm Lip}(X)$ precisely when $X$ is compact. In addition it was proved for
noncommutative Lipschitz algebras associated to noncommutative tori in
[15] and [16], and our proof here uses basically the same method. These results
also follow from an unpublished theorem of Rieffel which deals with the general
situation of a Lie group acting on a Banach space [13]. I do not know whether
that line of reasoning implies our current result (it is not obvious because
$\R^2$, the Lie group that appears here, is not compact).
\bigskip

\noindent {\bf Lemma 14.} {\it For any $\epsilon > 0$ there exists $N$ large
enough that $\|x - \sigma_n(x)\| \leq \epsilon$ for all $x \in
{\rm ball}(\L_\hbar)$ and $n \geq N$.}
\medskip

\noindent {\it Proof.} Recall that $\gamma_t =
\beta_{t'}^{-1}\alpha_{t'}^{-1}\beta_{t'}\alpha_{t'}$ where $t' =
\sqrt{t/c}$. Now for any $y \in {\rm ball}(\L_\hbar)$ we have
$$\|\alpha_{t'}(y) - y\|, \|\beta_{t'}(y) - y\| \leq t'\cdot\|d(y)\| \leq t',$$
so that $\|\gamma_t(x) - x)\| \leq 4t'$ for any $x \in {\rm ball}(\L_\hbar)$.
So we have
$$\|x - \sigma_n(x)\| \leq \tp\int_{-\pi}^\pi \|x - \gamma_t(x)\|K_n(t) dt
\leq \tp\int_{-\pi}^\pi 4\sqrt{t/c}K_n(t)dt,$$
and the last formula goes to zero as $n \to \infty$. This is what we needed
to show.\hfill\hal
\bigskip

\noindent {\bf Theorem 15.} {\it The unit ball of $\L_\hbar$ is compact in
operator norm.}
\medskip

\noindent {\it Proof.} Let $(x_k)$ be any sequence in ${\rm ball}(\L_\hbar)$;
we will find a convergent subsequence.
\medskip

As in Lemma 5, $Y_na_n(x_k)$ is multiplication by some function $g_n^k \in
L^\infty(\R\times\T\times\Z)$. Now
$$\eqalign{\|\alpha_r(Y_na_n(x_k)) - Y_na_n(x_k)\|
& = \|\alpha_r(a_n(x_k)) - a_n(x_k)\|\cr
&\leq
\tp\int_{-\pi}^\pi \|(\alpha_r(\gamma_t(x_k))-\gamma_t(x_k))e^{-int}\| dt\cr
& = \tp\int_{-\pi}^\pi \|\alpha_r(x_k) - x_k\| dt\cr
& \leq r\|\delta_1(x_k)\| \leq r.\cr}$$
Similarly $\|\beta_s(Y_na_n(x_k)) - Y_na_n(x_k)\| \leq s$, and
this implies that the function $g_n^k$ is Lipschitz with Lipschitz number at
most 1. Since $[0,2\pi]\times\T\times\{0\}$ is compact, we may choose
a subsequence $g_n^{k_j}$ which converges in sup norm on this set; by ($*$)
and ($\dag$) of Lemma 5
this implies that $g_n^{k_j}$ converges in sup norm on all of
$\R\times\T\times\Z$. Allowing $n$ to vary, finding successive subsequences 
for which $g_n^{k_j}$ converges, and diagonalizing, we get a subsequence
$(x_{j_k})$ of $(x_k)$ such that $(g_n^{k_j})$ converges in sup norm for all
$n$.
\medskip

Let $x$ be an weak operator cluster point of $(x_{j_k})$ and let $Y_n(a_n(x))$
be multiplication by $g_n$. Then $g_n$ is a cluster point of $(g_n^{k_j})$,
hence $g_n^{k_j} \to g_n$.
\medskip

Given $\epsilon > 0$, by Lemma 14 we can find a positive integer $N$ such
that $\|x - \sigma_N(x)\| \leq \epsilon$ and $\|x_{k_j} - \sigma_N(x_{k_j})\|
\leq \epsilon$ for all $j$. Taking $M$ large enough that $j \geq M$ implies
$$\|g_n - g_n^{k_j}\|_\infty \leq \epsilon/(2N + 1)$$
for all $|n| \leq N$, we get
$$\eqalign{\|\sigma_N(x) - \sigma_N(x_{k_j})\| &\leq
\sum_{n = -N}^N \Bigl(1 - {{|n|}\over{N+1}}\Bigr)\|a_n(x) - a_n(x_{k_j})\|\cr
&= \sum_{n= -N}^N \Bigl(1 - {{|n|}\over{N+1}}\Bigr)\|g_n - g_n^{k_j}\|\cr
&\leq \epsilon.\cr}$$
Thus
$$\|x - x_{k_j}\| \leq \|x - \sigma_N(x)\|
+ \|\sigma_N(x) - \sigma_N(x_{k_j})\| + \|\sigma_N(x_{k_j}) - x_{k_j}\|
\leq 3\epsilon$$
for $j \geq M$. So $x_{k_j} \to x$ in operator norm.\hfill\hal
\bigskip
\bigskip

\noindent {\bf 4. Further properties}
\bigskip

In this section we first identify the sub-Riemannian noncommutative Lipschitz
algebra $\L_\hbar$ defined in section 3 with a noncommutative H\"older
algebra. This generalizes the classical fact (see [10]) that the
sub-Riemannian metric on the Heisenberg manifolds is comparable to
the square root of the Riemannian metric in the $z$ direction. Then we use
work of Sauvageot [14] to establish the existence of a heat semigroup on
$N_\hbar$ and identify its generator with a noncommutative Laplacian.
\bigskip

\noindent {\bf Definition 16.} Let $A, B, C \in (0,1]$ and define
$\L_\hbar^{A,B,C}$ to be the set of $x \in N_\hbar$ for which there
exists a constant $K \geq 0$ such that
$$\|x - \alpha_r(x)\| \leq Kr^A,\qquad \|x - \beta_s(x)\| \leq Ks^B,
\qquad \|x - \gamma_t(x)\| \leq Kt^C$$
for all $r, s, t > 0$. Let $L(x) = L^{A,B,C}(x)$
be the smallest possible value of $K$ and norm $\L_\hbar^{A,B,C}$ by
$$\|x\|_{A,B,C} = {\rm max}(\|x\|, L(x)).\eqno{\hal}$$

Note that $A' \leq A, B' \leq B, C' \leq C$ implies $\L_\hbar^{A, B, C}
\subset \L_\hbar^{A', B', C'}$. Indeed, if $K = {\rm max}(L^{A,B,C}(x),2\|x\|)$
then
$$\|x - \alpha_r(x)\| \leq \cases{Kr^A \leq Kr^{A'}&for $r \leq 1$\cr
2\|x\| \leq Kr^{A'}& for $r \geq 1$,\cr}$$
and similarly for $\beta$ and $\gamma$, so that $L^{A',B',C'}(x) \leq
\max(L^{A,B,C}(x), 2\|x\|)$.
\medskip

Note also that $\L_\hbar^{1,1,1} = {\rm dom}(\alpha) \cap {\rm dom}(\beta)
\cap {\rm dom}(\gamma)$ by an obvious extension of Lemma 12. For this
reason we can realize $\L_\hbar^{1,1,1}$ as the domain of a derivation into
the Hilbert module $N_\hbar \oplus N_\hbar \oplus N_\hbar$ in the same way
that we treated $\L_\hbar$ in Definition 11. However, for non-unit values
of $A$, $B$, or $C$ Hilbert modules are not appropriate. We instead use a
construction from [16] to handle this case.
\bigskip

\noindent {\bf Definition 17.} Let
$$F = \bigoplus_{t > 0}^\infty (N_\hbar \oplus N_\hbar \oplus N_\hbar)$$
be the $l^\infty$ direct sum of von Neumann algebras. Regard it as a dual
operator $N_\hbar$-bimodule with left action given by the diagonal
embedding of $N_\hbar$ in $F$ and right action given by the embedding
$$x \mapsto \bigoplus_{t > 0}
(\alpha_t(x) \oplus \beta_t(x) \oplus \gamma_t(x)).$$
Define a map $d: \L_\hbar^{A,B,C} \to F$ by $d(x) = \bigoplus d_t(x)$ with
$$d_t(x) = {{x - \alpha_t(x)}\over{t^A}} \oplus {{x - \beta_t(x)}\over{t^B}}
\oplus {{x - \gamma_t(x)}\over{t^C}}.\eqno{\hal}$$

\noindent {\bf Theorem 18.} {\it The map $d$ and its domain $\L_\hbar^{A,B,C}$
have the following properties:
\medskip

{\narrower{
\noindent (a). $d$ is an unbounded derivation with weak*-closed graph.
\medskip

\noindent (b). $\L_\hbar^{A,B,C}$ is a dual Banach algebra. It
is densely contained in $D_\hbar$ and, if $C < 1$, it contains $S^c$.
\medskip

\noindent (c). The unit ball of $\L_\hbar^{A,B,C}$ is compact in operator
norm.
\medskip}}}

\noindent {\it Proof.} (a). It is routine to check that $d$ is a derivation.
To verify weak*-closure of the graph, let $(x_\lambda) \subset
\L_\hbar^{A,B,C}$ be a bounded net which weak operator converges to
$x \in N_\hbar$ and suppose $d(x_\lambda)$ weak operator converges to
$\bigoplus y_t^1 \oplus y_t^2 \oplus y_t^3 \in F$. Then restricting
attention to the $t$th summand of $F$, we have $d_t(x_\lambda) \to y_t^1
\oplus y_t^2 \oplus y_t^3$. Thus
$$(x_\lambda - \alpha_t(x_\lambda))/t^A \to y_t^1$$
(weak operator), but the left side also converges to $(x - \alpha_t(x))/t^A$.
An identical argument applies to the other two summands of $d_t(x)$ and so
we conclude that $d_t(x_\lambda) \to d_t(x)$. Boundedness of the net then
implies that $x \in \L_\hbar^{A, B, C}$ and $d(x_\lambda) \to d(x)$.
\medskip

(b). $\L_\hbar^{A, B, C}$ is a dual Banach algebra by the same easy argument
used in Theorem 13 (b). It is contained in $D_\hbar$ by Theorem 9.
Also, as remarked above it contains $L_\hbar^{1,1,1}$.
We claim that $\L_\hbar^{1,1,1}$ contains all $\Phi \in S^c$ which are zero
for all but a finite number of values of $p$. This implies density in $D_\hbar$
by Lemma 3.
\medskip

To prove the claim it suffices to consider only those $\Phi \in S^c$ for
which $\Phi(x,y,p) = 0$ unless $p = n$, for some fixed $n \in \Z$. For
such $\Phi$ the operator norm equals the $L^\infty$ norm, and so
$$\|\gamma_t(\Phi) - \Phi\| = \|(e^{int} - 1)\Phi\|_\infty \leq
nt\|\Phi\|_\infty.$$
Together with Lemma 8 this implies $\Phi \in \L_\hbar^{1,1,1}$, as claimed.
\medskip

Now suppose $C < 1$ and let $\Phi$ be any function in $S^c$. Choose
$N \geq 1 + 2/(1-C)$. Then by part (b) of Definition 1, there exists a constant
$K$ such that $|p^N\Phi(x,y,p)| \leq K$ for $x \in [0,2\pi]$ and all $y$ and
$p$. By part (a) of Definition 1, this implies that $p^N\Phi$ is bounded
on all of $\R\times\T\times\Z$. Now for any $t > 0$ define
$\Phi_q(x,y,p) = \delta_{p,q}\Phi(x,y,p)$ and
$\Phi' = \Phi - \sum_{|p| \leq t^{-1/(N-1)} + 1} \Phi_q$. 
Then $\|\Phi_q\|_\infty \leq q^{-N}K$ and so we can bound the
operator norm of $\Phi'$ by
$$\eqalign{\|\Phi'\|
&\leq \sum_{|q| > t^{-1/(N-1)} + 1} q^{-N}K\cr
& \leq 2\int_{t^{-1/(N-1)}}^\infty q^{-N}K\, dq\cr
& = 2{{(t^{-1/(N-1)})^{1-N}}\over{N - 1}}K\cr
& = 2tK/(N-1).\cr}$$
Thus $\|\gamma_t(\Phi') - \Phi'\| \leq 4tK/(N-1)$, so that
$\|\gamma_t(\Phi') - \Phi'\|/t^C$ is bounded for $t \leq 1$.
At the same time we have
$$|(\gamma_t(\Phi_q) - \Phi_q)(x,y,p)| = |(e^{iqt} - 1)\Phi_q(x,y,p)|
\leq |qt|\|\Phi\|_\infty \leq (t + t^{1 - 1/(N-1)})\|\Phi\|_\infty$$
for $|q| \leq t^{-1/(N-1)} + 1$, so that for $t \leq 1$ and $|q| \leq
t^{-1/(N-1)} + 1$ we have
$$\|\gamma_t(\Phi_q) - \Phi_q\| \leq 2t^{1-1/(N-1)}\|\Phi\|_\infty.$$
Hence
$$\|\gamma_t(\Phi - \Phi') - (\Phi - \Phi')\|
\leq (2t^{-1/(N-1)} + 3)\cdot 2t^{1-1/(N-1)}\|\Phi\|_\infty
\leq 12t^C\|\Phi\|_\infty.$$
We conclude that $\|\gamma_t(\Phi) - \Phi\|/t^C$ is bounded for
$t \leq 1$. But for $t \geq 1$ we have
$$\|\gamma_t(\Phi) - \Phi\| \leq 2\|\Phi\| \leq 2t^C\|\Phi\|,$$
so that $\|\gamma_t(\Phi) - \Phi\|/t^C$ is bounded for all $t > 0$. Thus
$\Phi \in \L_\hbar^{A,B,C}$.
\medskip

(c). This is proved in the same way as Theorem 15.\hfill\hal
\bigskip

\noindent {\bf Theorem 19.} {\it $\L_\hbar$ and $\L_\hbar^{1,1,1/2}$ are
identical as sets and have isomorphic norms.}
\medskip

\noindent {\it Proof.} $\L_\hbar^{1,1,1/2} \subset \L_\hbar$ is clear from
the definitions. Conversely, it was noted in the proof of Lemma 14 that
if $x \in {\rm ball}(\L_\hbar)$
then $\|\gamma_t(x) - x\| \leq 4\sqrt{t/c}$. Hence
$\|\gamma_t(x) - x\|/\sqrt{t}$ is bounded and we get $x \in
\L_\hbar^{1,1,1/2}$. So $\L_\hbar = \L_\hbar^{1,1,1/2}$ as sets.
\medskip

Isomorphism follows from the estimate
$$\|x\|_L \leq {\rm max}(\|x\|, \|\delta_1(x)\| + \|\delta_2(x)\|)
\leq {\rm max}(\|x\|, 2L(x)) \leq 2\|x\|_{A,B,C}$$
together with the open mapping theorem.\hfill\hal
\bigskip

As a final topic we consider a natural Laplacian and the heat semigroup it
generates. The following tool is needed.
\medskip

\noindent {\bf Definition 20.} Define $\tau: N_\hbar \to \C$ by
$$\tau(T) = \tp\int_0^1\int_0^{2\pi} \Phi(x,y,0)dxdy$$
where $\Phi$ is the function associated to $a_0(T)$ in Lemma 5.\hfill\hal
\bigskip

\noindent {\bf Theorem 21 ([12], p.\ 558).} {\it The map $\tau$ is a
faithful normal finite trace. It is invariant for the action $L$ given in
Definition 7.}\hfill\hal
\bigskip

(The fact that $\tau$ is normal follows from the fact that it is the
composition of two normal maps: integration of $\gamma_t(T)$ and integration
over $\T^2$.)
\bigskip

\noindent {\bf Proposition 22.} {\it The GNS representation of $N_\hbar$
associated to $\tau$ is unitarily equivalent to the restriction of its
original representation on $H = L^2(\R\times\T\times\Z)$ to $H' =
L^2([0,2\pi]\times\T\times\Z)$.}
\medskip

\noindent {\it Proof.} Define a map $\Psi: N_\hbar \to H'$ by
$\Psi(T) = T(\xi)$ where $\xi(x,y,p) = 1$ for $(x,y,p) \in
[0,2\pi]\times\T\times\{0\}$ and $\xi(x,y,p) = 0$ elsewhere. Then for any
$T \in N_\hbar$ we have $\tau(T) = \la T\xi, \xi\ra$, so that
$\Psi$ extends to a unitary map from the Hilbert space of the
GNS representation onto $H'$. It intertwines the action of $N_\hbar$ because
$$\Psi(TS) = TS(\xi) = T\Psi(S).\eqno{\hal}$$

Now we define a Laplacian on $S^c$, and apply a theorem of Sauvageot
to establish the existence of a heat semigroup.
\bigskip

\noindent {\bf Definition 23.} For $\Phi \in S^c$ define $\Delta\Phi$
by $\Delta\Phi = \delta_1(\delta_1(\Phi)) + \delta_2(\delta_2(\Phi))$.
Concretely, we have
$$\Delta\Phi = \Phi_{xx} -p^2c^2x^2\Phi(x,y,p) - 2ipcx\Phi_y + \Phi_{yy}.$$
The fact that $\Delta\Phi \in S^c$ follows immediately from the first
definition of $\Delta\Phi$,
or is a routine computation from the second.\hfill\hal
\bigskip

\noindent {\bf Theorem 24.} {\it The operators $e^{-t\Delta}$, $t \geq 0$,
form a weak-operator continuous semigroup of completely positive normal
contractions of $N_\hbar$.}
\medskip

\noindent {\it Proof.} First we show that $\delta_1 \oplus \delta_2$
is closed when regarded as an operator from $H'$ into $H' \oplus H'$.
\medskip

Identify $H'$ with the functions $\Phi$ on $\R\times\T\times\Z$ which satisfy
$$\Phi(x + k, y, p) = e^{ickpy}\Phi(x,y,p)$$
for all $k \in 2\pi\Z$ and whose restriction to $[0,2\pi]\times\T\times\Z$
is square integrable. Then for all $r, s, t \in \R$ the formula
$$(U_{(r,s,t)}\xi)(x,y,p) =
e^{ip(t + cs(x + \hbar p\mu - r))}\xi(x-r, y-s, p)$$
from Definition 7 defines a unitary operator on $H'$. Thus setting
$\alpha_r'(\xi) = U_{(r,0,0)}\xi$ and $\beta_s'(\xi) = U_{(0,s,0)}\xi$
gives us
two strongly continuous one-parameter unitary groups on $H'$, whose generators
$D_1$ and $D_2$ satisfy $D_i(\Phi) = \delta_i(\Phi)$, treating $\Phi$
respectively as an element of $H'$ and of $N_\hbar$.
\medskip

Since $D_1$ and $D_2$ are self-adjoint they are closed, hence
$$D = D_1\oplus D_2: H' \mapsto H'\oplus H'$$
is closed. Also $D^*(\xi\oplus 0) = D_1(\xi) = \delta_1(\xi)$ and
$D^*(0 \oplus \xi) = D_2(\xi) = \delta_2(\xi)$, so
$$D^*D(\Phi) = D^*(\delta_1(\Phi) \oplus \delta_2(\Phi))
= \delta_1(\delta_1(\Phi)) + \delta_2(\delta_2(\Phi)) = \Delta\Phi$$
for all $\Phi \in S^c$. The desired conclusion now follows from
([S], Corollary 3.5).\hfill\hal
\bigskip
\bigskip

[1] B.\ Abadie, ``Vector bundles'' over quantum Heisenberg manifolds, in
{\it Algebraic Methods in Operator Theory} (Birkh\"auser, 1994), 307-315.
\medskip

[2] {-------}, Isomorphism classes for irrational quantum Heisenberg
manifolds, preprint.
\medskip

[3] B.\ Abadie and R.\ Exel, Hilbert C*-bimodules over commutative C*-algebras
and an isomorphism condition for quantum Heisenberg manifolds, {\it Rev.\
Math.\ Phys.\ \bf 9} (1997), 411-423.
\medskip

[4] B.\ Abadie, S.\ Eilers, and R.\ Exel, Morita equivalence for
crossed products by Hilbert C*-bimodules, {\it Trans.\ Amer.\ Math.\ Soc.\ }
(to appear).
\medskip

[5] A.\ Bella\"iche, The tangent space in sub-Riemannian geometry, in [6],
1-78.
\medskip

[6] A.\ Bella\"iche and J.-J.\ Risler, eds., {\it Sub-Riemannian Geometry}
(Birkh\"auser, 1996).
\medskip

[7] O.\ Bratteli and D.\ W.\ Robinson, {\it Operator Algebras and Quantum
Statistical Mechanics, Vol.\ 1}, Springer-Verlag (1979).
\medskip

[8] W.\ L.\ Chow, \"Uber Systeme von Linearen partiellen
Differentialgleichungen erster Ordnung, {\it Math.\ Ann.\ \bf 117}
(1939), 98-105.
\medskip

[9] A.\ Connes, {\it Noncommutative Geometry} (Academic Press, 1994).
\medskip

[10] M.\ Gromov, Carnot-Carath\'eodory spaces seen from within, in [6], 79-323.
\medskip

[11] Y.\ Katznelson, {\it An Introduction to Harmonic Analysis}, Dover (1976).
\medskip

[12] M.\ A.\ Rieffel, Deformation quantization of Heisenberg manifolds,
{\it Comm.\ Math.\ Phys.\ \bf 122} (1989), 531-562.
\medskip

[13] {-------}, personal communication.
\medskip

[14] J.-L.\ Sauvageot, Quantum Dirichlet forms, differential calculus, and
semigroups, in {\it Quantum Probability and Applications V}, Springer LNM
\# 1442 (1990), 531-536.
\medskip

[15] N.\ Weaver, Lipschitz algebras and derivations of von Neumann algebras,
{\it J.\ Funct.\ Anal.\ \bf 139} (1996), 261-300.
\medskip

[16] {-------}, $\alpha$-Lipschitz algebras on the noncommutative torus,
to appear in {\it J.\ Operator Theory}.
\medskip

[17] {-------}, Operator spaces and noncommutative metrics, preprint.
\bigskip
\bigskip

\noindent Math Dept.

\noindent Washington University

\noindent St. Louis, MO 63130

\noindent nweaver@math.wustl.edu

\end